\documentclass[oneside,english,fullpage]{amsart}
\usepackage[T1]{fontenc}
\usepackage[latin9]{inputenc}
\usepackage{color}
\usepackage{amsthm}

\makeatletter
\numberwithin{equation}{section}
\numberwithin{figure}{section}
\theoremstyle{plain}
\newtheorem{thm}{\protect\theoremname}
  \theoremstyle{definition}
  \newtheorem{defn}[thm]{\protect\definitionname}
  \theoremstyle{plain}
  \newtheorem*{thm*}{\protect\theoremname}
  \theoremstyle{definition}
  \newtheorem*{example*}{\protect\examplename}

\usepackage{geometry}
\makeatother

\usepackage{babel}
  \providecommand{\definitionname}{Definition}
  \providecommand{\examplename}{Example}
  \providecommand{\theoremname}{Theorem}
\providecommand{\theoremname}{Theorem}

\begin{document}

\title[Condensation Method]{An Elementary Proof of Dodgson's Condensation Method for Calculating
Determinants}

\date{July 12, 2013}

\author{Mitch Main, Micah Donor, and R. Corban Harwood$^{1}$}

\address{$^{1}$Department of Mathematics and Applied Science, George Fox
University, Newberg, OR 97132, USA. Email: rharwood@georgefox.edu}
\begin{abstract}
In 1866, Charles Ludwidge Dodgson published a paper concerning a method
for evaluating determinants called the condensation method. His paper
documented a new method to calculate determinants that was based on
Jacobi\textquoteright s Theorem. The condensation method is presented
and proven here, and is demonstrated by a series of examples. The
condensation method can be applied to a number of situations, including
calculating eigenvalues, solving a system of linear equations, and
even determining the different energy levels of a molecular system.
The method is much more efficient than cofactor expansions, particularly
for large matrices; for a 5 x 5 matrix, the condensation method requires
about half as many calculations. Zeros appearing in the interior of
a matrix can cause problems, but a way around the issue can usually
be found. Overall, Dodgson\textquoteright s condensation method is
an interesting and simple way to find determinants. This paper presents
an elementary proof of Dodgson's method.
\end{abstract}

\maketitle
\global\long\def\R{\mathbb{R}}
\global\long\def\C{\mathbb{C}}
\global\long\def\Z{\mathbb{Z}}
\global\long\def\F{\mathbb{F}}

\section{Introduction}

Determinants are often used in linear algebra and other mathematical
situations. They first were developed in the 18th century, and have
since become ubiquitous in certain areas of mathematics. Although
small determinants are quite simple to calculate, the degree of difficulty
increases exponentially with the size of the matrix. An interesting
way of finding determinants that can be much simpler than traditional
methods came from an unlikely source. Reverend Charles Lutwidge Dodgson,
better known by his pseudonym Lewis Carroll as the creator of \textit{Alice
in Wonderland}, was also a recreational mathematician \cite{Article2 }.
He was especially fascinated by Euclid and published several works
on Euclidean geometry, although he did publish some other papers and
books on various topics, none of which became well known \cite{Article2 }.
In the year 1866, he published a paper in the Proceedings of the Royal
Society of London that first presented the condensation method to
the public \cite{Article2 }. In his paper Dodgson outlined his new
method for computing determinants, which he called ``far shorter
and simpler than any {[}method{]} hitherto employed''; he also proved
the method for 3 x 3 and 4 x 4 matrices \cite{Article1}. Relative
to using a cofactor expansion (the traditional method), the condensation
method can be quite a time-saver when dealing with larger matrices.
Dodgon's based his method on Jacobi's Theorem, a ``well-known theorem
in determinants'' \cite{Article1}. The method can be utilized to
improve calculation efficiency in many applications that involve finding
determinants, such as solving systems of equation using Cramer's Rule,
determining if a matrix is invertible, and finding the eigenvalues
of a matrix \cite{Article2 }\cite{Textbook1}. Determinants also
arise in other scientific fields; some examples are engineering, physics,
and as we shall see, chemistry. Dodgson's method is thus applicable
to a wide variety of situations and subjects. This paper will examine
the theoretical basis of the condensation method, as well as presenting
the method and some examples, and finally proving the method.

\section{Theory}

\subsection{Definitions}

The following definitions are relevant to the consideration of Charles
Dodgson's condensation method, and will be necessary to understand
the method and its applications.
\begin{defn}
For an $n\times n$ matrix, a \textit{minor} is any $(n-m)\times(n-m)$
matrix formed by deleting \textit{m} rows and \textit{m} columns from
the original matrix. \textcolor{black}{For} example, the italicized
entries in the matrix below form a $2\times2$ minor in the upper
left corner of the matrix.
\[
\left(\begin{array}{cccc}
{\color{red}\mathit{2}} & {\color{red}\mathit{1}} & -1 & -3\\
{\color{red}\mathit{1}} & {\color{red}\mathit{-2}} & 3 & 0\\
3 & 1 & 2 & -1\\
0 & -2 & 3 & 1
\end{array}\right)
\]

\end{defn}

\begin{defn}
A \textit{consecutive minor} is a minor in which the remaining rows
and columns were adjacent in the original matrix. Thus the example
above shows a consecutive minor, as does the example below.
\[
\left(\begin{array}{cccc}
{\normalcolor {\color{black}{\color{red}{\color{black}2}}}} & {\color{red}\mathit{1}} & {\color{red}\mathit{-1}} & -3\\
{\normalcolor 1} & {\color{red}\mathit{-2}} & {\color{red}\mathit{3}} & 0\\
3 & 1 & 2 & -1\\
0 & -2 & 3 & 1
\end{array}\right)
\]

\end{defn}

\begin{defn}
A \textit{complementary minor} is the $m\times m$ matrix diagonally
adjacent to the minor matrix. The italicized entries below form a
minor complementary to the original minor.
\[
\left(\begin{array}{cccc}
{\color{black}{\color{red}{\color{black}2}}} & {\color{black}{\color{black}{\color{red}{\color{black}1}}}} & -1 & -3\\
{\color{black}{\color{red}{\color{black}1}}} & {\color{red}\mathit{-2}} & {\color{red}\mathit{3}} & 0\\
3 & {\color{red}\mathit{1}} & {\color{red}\mathit{2}} & -1\\
0 & -2 & 1 & 1
\end{array}\right)
\]

\end{defn}

\begin{defn}
The \textit{interior} of a matrix is the matrix formed by deleting
the first and last rows and columns of the matrix. So for our example
matrix the interior is the italicized entries shown below.
\[
\left(\begin{array}{cccc}
{\color{black}{\color{red}{\color{black}2}}} & {\color{black}{\color{black}{\color{red}{\color{black}1}}}} & -1 & -3\\
{\color{black}{\color{red}{\color{black}1}}} & {\color{red}\mathit{-2}} & {\color{red}\mathit{3}} & 0\\
3 & {\color{red}\mathit{1}} & {\color{red}\mathit{2}} & -1\\
0 & -2 & 1 & 1
\end{array}\right)
\]

\end{defn}

\begin{defn}
A \textit{connected minor} is one in which all of the rows and columns
are adjacent. All of the examples shown have been connected minors.
\end{defn}

\begin{defn}
The \textit{adjugate} matrix is defined entrywise by $a'_{ij}=(-1)^{i+j}\cdot det[A_{ij}]$,
where $[A_{ij}]$ is the minor with row i and column j deleted. It
is denoted $A'$.
\end{defn}

\subsection{Jacobi's Theorem}

The foundation of Dodgson's method can be found in Jacobi's Theorem,
a result first stated by Jacobi in 1833, and then further developed
in papers published in 1835 and 1841, by which time the theorem was
fully realized \cite{Article3}. 
\begin{thm*}
Jacobi's Theorem. Let A be an $n\times n$ matrix, let $[A_{ij}]$
be an $m\times m$ minor of A, where $m<n$, let $[A'_{ij}]$ be the
corresponding $m\times m$ minor of A', and let $[A*_{ij}]$ be the
complementary $(n-m)\times(n-m)$ minor of A. Then:
\[
det[A'_{ij}]=det(A)^{m-1}\cdot det[A*_{ij}]
\]

\end{thm*}
Dodgson noted that by considering the case where \textit{m} is 2 and
dividing by $det[A*_{ij}]$ a useful algorithm for finding the determinant
of A could be obtained: 
\[
det(A)=\frac{det[A'_{ij}]}{det[A*_{ij}]}
\]

In their article about Charles Dodgson and his interesting method,
Rice and Torrence outline a proof of Jacobi's Theorem \cite{Article2 }.
Those interested can find it there.

\subsection{Condensation Method}

The condensation method developed by Charles Dodgson is as follows:
\begin{itemize}
\item First remove all zeros from the interior of A, using elementary row
and column operations. Call this matrix $A^{(0)}$.
\item Then take the determinants of the $2\times2$ consecutive minors to
form an $(n-1)\times(n-1)$ matrix $A^{(1)}$.
\item Next, take the determinants of the $2\times2$ consecutive minors
in $A^{(1)}$ to form an $(n-2)\times(n-2)$ matrix. Divide each term
by the corresponding entry from the interior of matrix $A^{(0)}$
to form $A^{(2)}$.
\item In general, given $A^{(k)}$, compute an $(n-k-1)\times(n-k-1)$ matrix
from the determinants of the $2\times2$ consecutive minors of $A^{(k)}$.
To produce $A^{(k+1)}$ divide each enty by the corresponding entry
in the interior of $A^{(k-1)}$. 
\item Repeat the previous step until a single number is obtained, which
is $det(A)$.
\end{itemize}

\section{Examples and Applications}

\subsection{Examples}

Now we shall consider some examples that illustrate the use of Dodgson's
condensation method.
\begin{example*}
Let 
\[
A=\begin{pmatrix}4 & 2 & 0 & -3\\
1 & 1 & 2 & 2\\
0 & -1 & 3 & -1\\
1 & 2 & 5 & 1
\end{pmatrix}
\]

Applying the condensation method, we obtain: 
\[
A^{(1)}=\left(\begin{array}{ccc}
\left|\begin{array}{cc}
4 & 2\\
1 & 1
\end{array}\right| & \left|\begin{array}{cc}
2 & 0\\
1 & 2
\end{array}\right| & \left|\begin{array}{cc}
0 & -3\\
2 & 2
\end{array}\right|\\
\left|\begin{array}{cc}
1 & 1\\
0 & -1
\end{array}\right| & \left|\begin{array}{cc}
1 & 2\\
-1 & 3
\end{array}\right| & \left|\begin{array}{cc}
2 & 2\\
3 & -1
\end{array}\right|\\
\left|\begin{array}{cc}
0 & -1\\
1 & 2
\end{array}\right| & \left|\begin{array}{cc}
-1 & 3\\
2 & 5
\end{array}\right| & \left|\begin{array}{cc}
3 & -1\\
5 & 1
\end{array}\right|
\end{array}\right)=\begin{pmatrix}2 & 4 & 6\\
-1 & 5 & -8\\
1 & -11 & 8
\end{pmatrix}
\]

So then 
\[
A^{(2)}*=\begin{pmatrix}14 & -62\\
6 & -48
\end{pmatrix}
\]

Dividing each entry by the corresponding entry of the interior of
A - in other words, taking the upper left entry of $A^{(2)}*$, 14,
and dividing by the upper left entry of the interior of $A$, 1, then
continuing in this fashion, we get
\[
A^{(2)}=\begin{pmatrix}14 & -31\\
-6 & -16
\end{pmatrix}
\]

Thus, $A^{(3)}*=(-410)$, and dividing by the interior of $A^{(1)}$
(which is 5), we get the result $A^{(3)}=-82=det(A)$.
\end{example*}

\begin{example*}
Let us consider another example. Let 
\[
A=\begin{pmatrix}0 & 1 & 0 & 4\\
-1 & 3 & 6 & -3\\
5 & 1 & 2 & 0\\
-2 & 1 & -1 & 1
\end{pmatrix}
\]

Then we find 
\[
A^{(1)}=\begin{pmatrix}1 & 6 & -24\\
-16 & \underline{{\color{black}0}} & 6\\
7 & -3 & 2
\end{pmatrix}
\]

This matrix has a zero in the interior, which will cause problems
later on. We can use row operations to change A so we are able to
find the determinant. We rearrange A by placing the first row at the
bottom and moving all the other rows up one, for a total of 3 row
interchanges. We now have the matrix 
\[
B=\begin{pmatrix}-1 & 3 & 6 & -3\\
5 & 1 & 2 & 0\\
-2 & 1 & -1 & 1\\
0 & 1 & 0 & 4
\end{pmatrix}
\]

Now after we take determinants of consecutive minors, we get
\[
B^{(1)}=\left(\begin{array}{ccc}
-16 & 0 & 6\\
7 & -3 & 2\\
-2 & 1 & -4
\end{array}\right)
\]

As you can see, the zero is no longer in the interior of the matrix,
so we can proceed in our calculations. Continuing on , we eventually
obtain $det(B)=163$. Since we initially performed 3 row interchanges,
we must multiply our answer by $(-1)^{3}$, and we find that $det(A)=-163$.
\end{example*}

\subsection{Applications}

The condensation method can be applied to many situations where calculation
of a determinant is needed. One such situation is in quantum mechanics,
where determinants can be used to determine the energy levels of molecules
(specifically $\pi$ electron systems) using Molecular Orbital Theory.
Using some experimental data, as well as theory developed for simple
systems, we can generate what is called a system of secular equations
that can be solved by setting the determinant of the system equal
to zero. The solutions to this system give the energy levels of the
particular molecule being studied; knowing the energy levels is crucial
to understanding the bonding and other electronic characteristics
of the molecule.The energy levels themselves represent the energy
values that electrons can have in the system in question. Since energy
is quantized at the molecular level, only the values given by the
solutions to the secular determinant are permitted. For example, in
linear $H_{3}$ the secular determinant using the Huckel approximations
is 
\[
\left|\begin{array}{ccc}
\alpha-E & \beta & 0\\
\beta & \alpha-E & \beta\\
0 & \beta & \alpha-E
\end{array}\right|=0
\]

In this determinant, $\alpha$ is called a Coulomb integral, and it
represents the energy of the electron when it is held by one of the
atoms; it has a negative value. The parameter $\beta$ is called a
resonance integral, and it represents the electron being between two
atoms, i.e. the wavefunctions of the two atoms overlapping. It is
also negative. The rows and columns of the determinant denote the
atoms in the molecule, so entry 1-1 is atom 1 combining with itself,
entry 1-2 is atoms 1 and 2 overlapping, etc. Entries 1-3 and 3-1 are
zero, since atoms 1 and 3 are not adjacent and thus cannot have their
wavefunctions overlap. We can use the condensation method to find
this determinant: 
\[
A^{(1)}=\left(\begin{array}{cc}
(\alpha-E)^{2}-\beta^{2} & \beta^{2}\\
\beta^{2} & (\alpha-E)^{2}-\beta^{2}
\end{array}\right)
\]

We then get
\[
A^{(2)}*=(\alpha-E)^{4}-2\beta^{2}(\alpha-E)^{2}
\]

So then 
\[
A^{(2)}=(\alpha-E)^{3}-2\beta^{2}(\alpha-E)=det(A)
\]

Setting this equal to zero we can then solve for the energy levels,
obtaining
\[
E=\alpha,\ \alpha\pm\sqrt{2}\beta
\]

Thus we know what the three energy levels are for linear $H_{3}$,
and we can use this data to find the bonding energy of $H_{3}$, which
can give us information about the stability and reactivity of the
molecule. For a more in-depth discussion of this topic, refer to Atkins
and De Paula's \textit{Physical Chemistry} \cite{Textbook3}.

\section{Proof of the Condensation Method}

In his book on the Alternating Sign Matrix Conjecture, David Bressoud
outlines the general idea behind proving the condensation method,
but does not go into detail as to the actual mechanics of the proof
\cite{Textbook2}. A combinatorial proof has also been published \cite{Article4}.
However, it is rather complicated, and so we endeavored to discover
a clear proof of Charles Dodgson's method.
\begin{proof}
We shall prove the condensation method using mathematical induction. 

In all cases, assume that there are no zeros in the interior of the
matrix, or that they have been removed using row operations prior
to beginning the process of condensation.

We will notate $A_{ij}^{kl}$ as the $2\times2$ minor corresponding
to the \textit{i-}th and \textit{j}-th rows and \textit{k}-th and
\textit{l}-th columns of $A$, and $a'_{ij}$ as the entry of the
adjugate matrix in the \textit{i}-th row and \textit{j}-th column.
In general, given $A_{m}^{n}$, $m$ represents the rows and $n$
the columns.

\medskip{}

\textbf{Base Case}: $3\times3$ matrix. 

Consider an arbitrary $3\times3$ matrix:
\[
A=\left(\begin{array}{ccc}
a_{11} & a_{12} & a_{13}\\
a_{12} & a_{22} & a_{23}\\
a_{13} & a_{23} & a_{33}
\end{array}\right)
\]

Using the condensation method, after simplification we get
\[
det(A)=\frac{det(A_{12}^{12})\cdot det(A_{23}^{23})-det(A_{12}^{23})\cdot det(A_{23}^{12})}{a_{22}}
\]

Applying Jacobi's Theorem, after simplifying we get
\[
det(A)=\frac{a'_{11}a'_{33}-a'_{13}a'_{31}}{a_{22}}=\frac{det(A_{12}^{12})\cdot det(A_{23}^{23})-det(A_{12}^{23})\cdot det(A_{23}^{12})}{a_{22}}
\]

Thus we can see that using the condensation method is equivalent to
using Jacobi's Theorem for a $3\times3$ matrix, and thus by Jacobi's
Theorem the condensation method is valid for $3\times3$ matrices.
\medskip{}

\textbf{Induction Hypothesis}: assume that the condensation method
is valid for \textit{$k\times k$} matrices ($k-1$ steps).

\medskip{}

It is vital to note that when using the condensation method, after
1 step we have a new matrix composed of the determinants of the $2\times2$
connected minors, after 2 steps we have the $3\times3$ determinants,
and so forth; after $i$ steps we have a matrix composed of the determinants
of the $(i+1)\times(i+1)$ connected minors. 

Now consider a $(k+1)\times(k+1)$ matrix. By our induction hypothesis,
we can condense the matrix to a $2\times2$ matrix in $k-1$ steps.
We are now left with a $2\times2$ matrix, each entry of which is
the determinant of one of the four $k\times k$ connected minors of
the original matrix:
\[
\left(\begin{array}{cc}
det(A_{1\leq i\leq k}^{1\leq j\leq k}) & det(A_{1\leq i\leq k}^{2\leq j\leq k+1})\\
det(A_{2\leq i\leq k+1}^{1\leq j\leq k}) & det(A_{2\leq i\leq k+1}^{2\leq j\leq k+1})
\end{array}\right)
\]

We can now complete the last step of the condensation method: take
the determinant of this $2\times2$ matrix, then divide the result
by the interior of the previous $3\times3$ matrix, which will be
the determinant of the interior of the original matrix. The result
is thus:
\begin{eqnarray*}
\frac{det(A_{1\leq i\leq k}^{1\leq j\leq k})\cdot det(A_{2\leq i\leq k+1}^{2\leq j\leq k+1})-det(A_{2\leq i\leq k+1}^{1\leq j\leq k})\cdot det(A_{1\leq i\leq k}^{2\leq j\leq k+1})}{det(A_{2\leq i\leq k}^{2\leq j\leq k})} & = & \frac{a'_{11}\cdot a'_{(k+1)(k+1)}-a'_{1(k+1)}\cdot a'_{(k+1)1}}{det(A*_{1(k+1)})}\\
 & = & \frac{det(A'_{1(k+1)})}{det(A*_{1(k+1)})}\\
 & = & det(A)\ (Jacobi's\ Thm.)
\end{eqnarray*}

Thus, by Jacobi's Theorem, the condensation method is valid for all
$n\times n$ matrices. 
\end{proof}

\section{Discussion}

The condensation method, first discovered by Charles Dodgson in 1866,
provides a more efficient way of calculating determinants of large
matrices. For a 5 x 5 matrix, it reduces the amount of computations
by nearly half compared to a standard cofactor expansion method. Dodgson
based his method on Jacobi's Theorem, which after some rearrangement
provided a useful algorithm for finding determinants.The main weakness
of the condensation method is zeros in the interior of the matrix,
which create issues with the necessary divisions. However, that hurdle
can often be circumvented by transforming the matrix using elementary
row operations, and in cases where there are too many zeros to use
the condensation method a traditional cofactor expansion is often
very easy to use. The condensation method can be applied in a variety
of situations, such as solving a system of linear equations with Cramer's
Rule, finding the eigenvalues of a matrix, determining invertibility,
and so forth. It can also be used in a quantum mechanical setting
to find the energy levels of a $\pi$ electron system, as was demonstrated
above. The mtheod preserves factorizations somewhat better than the
cofactor method, which is important in certain situations, particularly
finding the eigenvalues with the characteristic polynomial. Dodgson's
method is not well known due to his general mathematical obscurity,
but it can be extremely useful for finding determinants, especially
of larger matrices.

\end{document}